\documentclass{elsart}

\usepackage{amssymb,latexsym,amsmath,amsfonts}

\newtheorem{theorem}{Theorem}[section]

\newtheorem{corollary}[theorem]{Corollary}

\newtheorem{proposition}[theorem]{Proposition}

\begin{document}

\begin{frontmatter}



\title{On singular value inequalities for matrix means}


\author[1]{Raluca Dumitru\corauthref{Raluca Dumitru}}
\address[1]{Department of Mathematics and Statistics, University of North Florida; Institute of Mathematics of the Romanian Academy, Bucharest, Romania}
\ead{raluca.dumitru@unf.edu}
\corauth[Raluca Dumitru]{Department of Mathematics and Statistics, University of North Florida, 1 UNF Dr., Jacksonville, FL 32224; Telephone: 9046203740; Fax: 9046202818}

\author[2]{Rachel Levanger}
\address[2]{Department of Mathematics , Rutgers University}
\ead{rachel.levanger@gmail.com}

\author[3]{Bogdan Visinescu}
\address[3]{Department of Mathematics and Statistics, University of North Florida; Institute of Mathematics of the Romanian Academy, Bucharest, Romania}
\ead{b.visinescu@unf.edu}

\begin{abstract}
For positive semidefinite $n\times n$ matrices $A$ and $B$, the singular value inequality $(2+t)s_{j}(A^{r}B^{2-r}+A^{2-r}B^{r})\leq 2s_{j}(A^{2}+tAB+B^{2})$ is shown to hold for $r=\frac{1}{2}, 1, \frac{3}{2}$ and all $-2<t\leq 2$.
\end{abstract}


\end{frontmatter}



\newpage

\section{Introduction}

The arithmetic-geometric mean inequality for positive real numbers $a$ and $b$, $ab\leq (a^{2}+b^{2})/2$ has been proved by Bhatia and Kittaneh \cite{bk} to hold for singular values of arbitrary $n \times n$ matrices $A$ and $B$:
\begin{equation}\label{ag}
2s_{j}(AB^{\ast})\leq s_{j}(A^{\ast}A+B^{\ast}B)
\end{equation}
for all $j = 1, 2, ..., n$.
Using an additional $n \times n$ matrix $X$, Bhatia and Davis \cite{bd} obtained the following operator norm inequality:
$$2\|AXB^{\ast}\|\leq \| A^{\ast}AX+XB^{\ast}B\|$$
for all unitarily invariant norms, an inequality which was generalized by \\ X. Zhan \cite{zhan} as stated below.

Let $A,B,X\in M_{n}(C)$ with $A$ and $B$ positive semidefinite, $1\leq 2r\leq 3$, and $-2<t\leq 2$. Then for any unitarily invariant norm,
$$(2+t)\|A^{r}XB^{2-r}+A^{2-r}XB^{r}\|\leq 2\|A^{2}X+tAXB+XB^{2}\|.$$

X. Zhan conjectured \cite{zhanconj} that the above inequality has a similar form for singular values. Namely, for $X=I_{n}$ above, and for $A,B\geq 0$, $1\leq 2r\leq 3$, $-2<t\leq 2$ and all $j=1,2,...,n$, then
\begin{equation}\label{conj}
(2+t)s_{j}(A^{r}B^{2-r}+A^{2-r}B^{r})\leq 2s_{j}(A^{2}+tAB+B^{2}).
\end{equation}

These inequalities have been proved to hold for $t=0$ and all $0\leq r\leq 2$, a result due to Audenaert \cite{auden}. In this paper, we prove the inequalities in equation (\ref{conj}) for $r=\frac{1}{2},1,\frac{3}{2}$ and all $-2<t\leq 2$.

\section{Preliminaries}\label{prelim}

With $A$ a positive semidefinite matrix (we use the standard notation $A\geq 0$), denote by $\lambda_{j}(A)$ and $s_{j}(A)$ its eigenvalues and singular values, respectively, arranged in non-increasing order. Denote by $u_{j}(A)$ the eigenvectors of A corresponding to the eigenvalues $\lambda_{j}(A)$. 

Let $M_{0}$ and $M_{1}$ be two positive semidefinite $n\times n$ matrices. Let $M(t)=M_{0}+tM_{1}$, where $t$ is a scalar parameter. Then ( \cite{kato}, Chapter 2, Section 1) the eigenvalues $\lambda_{j}(M(t))$, $j=1,...,s$, of $M(t)$ are branches of analytic functions of $t$ with only algebraic singularities and the number $s$ of eigenvalues of $M(t)$ is independent of $t$, with the exception of a finite number of points $t$, called exceptional points. More precisely, there are exactly two possibilities. If $s=n$ (when these analytic functions are all distinct), then $M(t)$ has simple spectrum for all non-exceptional points $t$. If, on the other hand, some of these analytic functions are identical, then $s<n$, and $M(t)$ is called permanently degenerate. In conclusion, $\lambda_{j}(M(t))$ are everywhere continuous functions of $t$ and they are differentiable everywhere, except maybe at a finite number of points (the exceptional points). 

Let us consider next the derivative with respect to $t$ of these functions, whenever the derivative exists.

In what follows, we will use the notations $u_{j}(t)=u_{j}(M(t))$ and $\lambda_{j}(t)=\lambda_{j}(M(t))$.For simplicity, we will look at the derivative of $\lambda_{j}(M(t))$ at $t=0$. We will consider three cases.

Case 1. We assume that the eigenvalues of $M_{0}$ are all simple. 

Then the eigenvalues of $M(t)$ are also simple for small enough values of t, say $t\in(-a,a)$, for some $a$. This follows from Weyl's inequalities,

$$\lambda_{j}(M_{0})+\lambda_{n}(tM_{1})\leq \lambda_{j}(M_{0}+tM_{1})\leq \lambda_{j}(M_{0})+\lambda_{1}(tM_{1})$$

Hence, if $t\|M_{1}\|$ is small enough, namely strictly less than one half the minimum distance between all pairs of eigenvalues of $M_{0}$, then the eigenvalues of $M_{0}+tM_{1}$ will be simple.

Therefore, on the interval $t\in(-a,a)$, every eigenvalue $\lambda_{j}(t)$ has a unique eigenvector $u_{j}(t)$, up to a constant multiple. Furthermore, since the zeros of the characteristic polynomial of $M(t)$ are simple (hence the polynomial has nonzero derivative at these zeros) then by applying the implicit function theorem, we get that the eigenvalues $\lambda_{j}(t)$ are smooth in $(-a,a)$. Since the eigenvectors $u_{j}(t)$ are determined up to a scalar by the equations $(M(t)- \lambda_{j}(t)I_{n})u_{j}(t)=0$ and $u_{j}(t)^{\ast}u_{j}(t)=1$, the implicit function theorem allows us to locally select the functions $u_{j}(t)$ to also be smooth.

Let now $t\in(-a,a)$. By differentiating the equation
$$M(t)u_{j}(t)=\lambda_{j}(t)u_{j}(t)$$
we obtain the relation
\begin{equation}\label{[1]}
M'(t)u_{j}(t)+M(t)u'_{j}(t)=\lambda'_{j}(t)u_{j}(t)+\lambda_{j}(t)u'_{j}(t)
\end{equation}
and, by taking the inner product with $u_{j}(t)$ in equation (\ref{[1]}), we obtain
\begin{equation}\notag
u_{j}^{\ast}(t)M'(t)u_{j}(t)+u_{j}^{\ast}(t)M(t)u'_{j}(t)=\lambda'_{j}(t)u_{j}^{\ast}(t)u_{j}(t)+\lambda_{j}(t)u_{j}^{\ast}(t)u'_{j}(t).
\end{equation}

Since $M(t)$ is Hermitian, we have $u_{j}^{\ast}(t)M(t)=\lambda_{j}(t)u_{j}^{\ast}(t)$ and, using $u_{j}^{\ast}(t)u_{j}(t)=1$, we get
\begin{equation}\label{[2]}
\frac{d}{dt}\lambda_{j}(t)=u_{j}^{\ast}(t)M'(t)u_{j}(t).
\end{equation}

In particular, the derivative at 0 of the eigenvalue function in the case when $M_{0}$ has simple eigenvalues is given by
\begin{equation}
\frac{d}{dt}\bigg| _{t=0}\lambda_{j}(t)=u_{j}^{\ast}(0)M_{1}u_{j}(0).
\end{equation}

We assume next that $M_{0}$ has degenerate eigenvalues. Since the degeneracy of the eigenvalues of $M(t)$ is either accidental (for isolated values of $t$, such as $t=0$ here) or permanent (for all values of $t$), there are two cases to consider.

We will first consider the case when $M_{1}$ is such that it removes the degeneracy of $M_{0}$ for small enough values of $t$, so $M(t)$ has an exceptional point at $t=0$. The second case will be when $M(t)$ is permanently degenerate.

The problem is, in both these two cases, that the eigenvectors $u_{j}(t)$ are no longer unique.

Case 2. Assume that $t=0$ is an exceptional point. 

Then for $t$ small enough, $t\in(-a,a)\backslash {0}$, the eigenvalues of $M(t)$ are simple, and therefore the corresponding eigenvectors are unique. Hence $\lambda_{j}(t)$ is differentiable for all values of $t$ in $(-a,a)\backslash {0}$ and equation \ref{[2]} above still holds. Note that $\lambda_{j}(t)$ might not be differentiable at $t=0$, but it is continuous.

Case 3. Assume that $M(t)$ is permanently degenerate. 

Let now $a$ be the largest positive value for which $M(t)$ has no exceptional point in the interval $(-a,a)$. Let $\lambda_{j}(t)$ be an eigenvalue function such that $\lambda_{j}(0)$ has multiplicity $m$. 

Using  \cite{kato}(Chapter 2, equations 2.3, 2.21 and 2.34), we get

\begin{equation}\label{deriv}
\frac{d}{dt}\lambda_{j}(t)=\sum_{n=1}^{\infty}nt^{n-1}\frac{1}{mn}Tr(M'(0)P_{j}^{(n-1)}),
\end{equation}
where $M'(0)=\frac{d}{dt}\bigg| _{t=0}M(t)=M_{1}$ and $P_{j}(t)=\sum_{n=0}^{\infty}t^{n}P_{j}^{(n)}$ denotes the projection onto the eigenspace generated by $\lambda_{j}(t)$.

Note that, in this setting, there is no splitting of $\lambda$ in the interval $(-a,a)$, so that the $\lambda$-group consists of a single eigenvalue of multiplicity $m$. Hence the weighted mean of the $\lambda$-group, $\widehat{\lambda}_{j}(t)$ (used in equations 2.21 and 2.34) is the same as $\lambda_{j}(t)$.

Therefore,
\begin{align}\label{deriv}
\frac{d}{dt}\lambda_{j}(t)&=\frac{1}{m}\sum_{n=1}^{\infty}Tr(M'(0)t^{n-1}P_{j}^{(n-1)})\notag\\ 
			  &=\frac{1}{m}Tr(M'(0)P_{j}(t)),
\end{align}
for all t in the interval $(-a,a)$.

\section{The Behavior of an Eigenvalue Function}

\begin{theorem}\label{eig} For $A,B\in M_{n}(C)$, $A,B\geq 0$, $j=1,2,...,n$ and $t\in (-2,\infty)$, the function
$$f(t)=\frac{1}{2+t}\lambda_{j}(A^{2}+B^{2}+\frac{t}{2}AB+\frac{t}{2}BA)$$
is non-increasing.
\end{theorem}
\begin{pf}
Let $A,B\geq 0$. Note that $f(t)$ is continuous everywhere on $(-2,\infty)$ and it is differentiable everywhere except maybe at a finite number of points (the exceptional points). Let $t_{0}\in(-2,\infty)$ . We prove next that there exists an interval centered at $t_{0}$ where the function $f$ is non-increasing.

Let $M(t)=A^{2}+B^{2}+\frac{t}{2}AB+\frac{t}{2}BA$. We will consider three cases, according to whether $M(t)$ is (permanently) degenerate or not.

Case1. Assume that all eigenvalues of $M(t_{0})$ are simple. Then the eigenvalues of $M(t)$ are simple for values of t in a small enough neighborhood of $t_{0}$, say $t\in (t_{0}-a, t_{0}+a)$.

Then $f$ is differentiable everywhere on $(t_{0}-a, t_{0}+a)$ and using equation \ref{[2]} above, we obtain 

\begin{align}\label{deriv1}
f'(t)&=\frac{\frac{d}{dt}\lambda_{j}(t)(2+t)-\lambda_{j}(t)}{(2+t)^{2}}\\ \notag
                  &=\frac{u_{j}^{\ast}(t)\frac{AB+BA}{2}u_{j}(t)(2+t)-\lambda_{j}(t)}{(2+t)^{2}}\\ \notag
                  &=\frac{u_{j}^{\ast}(t)(\frac{AB+BA}{2}(2+t)-M(t))u_{j}(t)}{(2+t)^{2}}\\ \notag
                  &=-\frac{u_{j}^{\ast}(t)(A-B)^{2}u_{j}(t)}{(2+t)^{2}},\notag
\end{align}
for all $t\in (t_{0}-a, t_{0}+a)$.

Since $(A-B)^{2}\geq 0$, then for all $u\in C^{n}$ we have $\langle(A-B)^{2}u,u\rangle \geq 0$. 

In particular, this implies that $u_{j}^{\ast}(t)(A-B)^{2}u_{j}(t)\geq 0$, and therefore $f'(t)\leq 0$ , so $f$ is non-increasing on $(t_{0}-a, t_{0}+a)$.

Case 2. Assume that $t_{0}$ is an exceptional point for $M(t)$, hence $f$ might not be differentiable at $t_{0}$. 

There is however a small interval centered at $t_{0}$, say $(t_{0}-a, t_{0}+a)$, such that $f$ is differentiable everywhere on $(t_{0}-a, t_{0}+a) \backslash \{t_{0}\}$. Then the derivative of $f$ will be computed as in equation \ref{deriv1} and hence $f'(t)\leq 0$ on $(t_{0}-a, t_{0}+a)\backslash \{t_{0}\}$. Since $f$ is continuous everywhere, we conclude that $f$ is non-increasing on $(t_{0}-a, t_{0}+a)$.

Case 3. Assume that $M(t)$ is permanently degenerate and let $a$ be the largest positive value for which $M(t)$ has no exceptional point in the interval $(t_{0}-a,t_{0}+a)$. Let $\lambda_{j}(t)$ be an eigenvalue function such that $\lambda_{j}(t_{0})$ has multiplicity $m$.

Then $f(t)$ is everywhere differentiable on $(t_{0}-a, t_{0}+a)$, however its derivative cannot be computed in the same way as in equation \ref{deriv1} since the corresponding eigenvectors are not unique anymore.

Using equation \ref{deriv}, we obtain
\begin{align}
f'(t)&=\frac{(2+t)\frac{d}{dt}\lambda_{j}(t)-\lambda_{j}(t)}{(2+t)^{2}}\notag\\
     &=\frac{(2+t)\frac{1}{m}Tr(M'(0)P_{j}(t))-\frac{1}{m}Tr(M(t)P_{j}(t))}{(2+t)^{2}}\notag\\
     &=\frac{\frac{1}{m}Tr(((2+t)M'(0)-M(t))P_{j}(t))}{(2+t)^{2}}\notag\\
     &=\frac{-\frac{1}{m}Tr((A-B)^{2}P_{j}(t))}{(2+t)^{2}}.\notag
\end{align}
Hence $f'(t)\leq 0$ on $(t_{0}-a, t_{0}+a)$ and we conclude again that $f$ is non-increasing on $(t_{0}-a, t_{0}+a)$.
\qed
\end{pf}

Using Theorem \ref{eig} for  $t \in (-2,2]$, we obtain in particular the following corollary.

\begin{corollary}\label{cor} For $A,B\in M_{n}(C)$, $A,B\geq 0$, $j=1,2,...,n$ and $t\in (-2,2]$, we have
\begin{enumerate}
\item $A^{2}+B^{2}+\frac{t}{2}AB+\frac{t}{2}BA\geq 0$ and
\item $\frac{1}{2+t}s_{j}(A^{2}+B^{2}+\frac{t}{2}AB+\frac{t}{2}BA)\geq \frac{1}{4}s_{j}(A+B)^{2}.$
\end{enumerate}
\end{corollary}

\section{On Zhan's Conjecture}

We prove first that X. Zhan's conjecture (equation (\ref{conj})) holds for $r=\frac{1}{2},\frac{3}{2}$ and all $-2<t\leq 2$.

\begin{proposition}\label{prop}
For $A,B\in M_{n}(C)$, $A,B\geq 0$, $j=1,2,...,n$ and $t\in (-2,2]$, we have
$$(2+t)s_{j}(A^{\frac{1}{2}}B^{\frac{3}{2}}+A^{\frac{3}{2}}B^{\frac{1}{2}})\leq 2s_{j}(A^{2}+tAB+B^{2}).$$
\end{proposition}

\begin{pf}
Since $A^{2}+B^{2}+\frac{t}{2}AB+\frac{t}{2}BA=Re(A^{2}+tAB+B^{2})$, by using \cite{bhatia} Proposition III.5.1 we get
$$\lambda_{j}(A^{2}+B^{2}+\frac{t}{2}AB+\frac{t}{2}BA)\leq s_{j}(A^{2}+tAB+B^{2})$$
which, by Corollary \ref{cor} (1), is the same as 
\begin{equation}\label{ineq1}
s_{j}(A^{2}+B^{2}+\frac{t}{2}AB+\frac{t}{2}BA)\leq s_{j}(A^{2}+tAB+B^{2}).
\end{equation}
Using Corollary \ref{cor} (2), we conclude that 
\begin{equation}\label{ineq2}
\frac{2+t}{4}s_{j}((A+B)^{2})\leq s_{j}(A^{2}+tAB+B^{2}).
\end{equation}
Bhatia and Kittaneh \cite{bk1} proved that 
$$2s_{j}(A^{\frac{1}{2}}B^{\frac{3}{2}}+A^{\frac{3}{2}}B^{\frac{1}{2}})\leq  s_{j}((A+B)^{2}),$$
which, combined with equation (\ref{ineq2}), proves the desired result.
\qed
\end{pf}

Our next result shows that Zhan's conjecture (equation (\ref{conj})) holds for $r=1$ and all $-2<t\leq 2$.

\begin{proposition}
For $A,B\in M_{n}(C)$, $A,B\geq 0$, $j=1,2,...,n$ and $t\in (-2,2]$, we have
$$(2+t)s_{j}(AB)\leq s_{j}(A^{2}+tAB+B^{2}).$$
\end{proposition}

\begin{pf}
Note that $4s_{j}(AB)\leq s_{j}(A+B)^{2}$ , an inequality recently proved by Drury in \cite{drury}. Using this inequality together with Corollary \ref{cor} (2) and Equation (\ref{ineq1}) we obtain

\begin{align}\notag
(2+t)s_{j}(AB)&\leq \frac{2+t}{4}s_{j}(A+B)^{2}\\ \notag
              &\leq s_{j}(A^{2}+B^{2}+\frac{t}{2}AB+\frac{t}{2}BA)\\ \notag
              &\leq s_{j}(A^{2}+tAB+B^{2}).\notag
\end{align}
\qed
\end{pf}


\begin{thebibliography}{00}





\bibitem{auden}K. Audenaert, \emph{A singular value inequality for Heinz means}, Linear Algebra Appl., 422(2007), 279-283.
\bibitem{bhatia}R. Bhatia, \emph{Matrix Analysis}, Springer, 1996.
\bibitem{bd}R. Bhatia, C. Davis, \emph{More matrix forms of the arithmetic-geometric mean inequality}, SIAM J. Matrix Anal. Appl., 14(1993), 132-136.
\bibitem{bk1}R. Bhatia, F. Kittaneh, \emph{Notes on matrix arithmetic-geometric mean inequality}, Linear Algebra Appl., 308(2000), 203-211.
\bibitem{bk}R. Bhatia, F. Kittaneh, \emph{On the singular values of a product of operators}, SIAM J. Matrix Anal. Appl., 11(1990), 272-277.
\bibitem{drury}S. W. Drury, \emph{On a question of Bhatia and Kittaneh}, Linear Algebra Appl., 437(2012), 1955-1960.
\bibitem{kato}T. Kato, \emph{Perturbation Theory for linear operators},  Springer-Verlag, Berlin, Heidelberg, New York, 1980.
\bibitem{zhan}X. Zhan, \emph{Inequalities for unitarily invariant norms}, SIAM J. Matrix Anal. Appl., 20(1999), 466-470.
\bibitem{zhanconj}X. Zhan, \emph{Some research problems on the Hadamard product and singular values of matrices}, Linear and Multilinear Algebra, Vol. 47, 191-194.




\end{thebibliography}
\end{document}